\documentclass{gtart}
\usepackage{amsmath,amssymb,amscd,enumerate}

\input gtoutput
\volumenumber{3}\papernumber{1}\volumeyear{1999}
\pagenumbers{1}{20}\published{15 March 1999}
\proposed{Frances Kirwan}\seconded{Simon Donaldson, Robion Kirby}
\received{19 May 1998}\revised{27 November 1998}
\accepted{16 February 1999}

\newtheorem{thm}{Theorem}
\newtheorem*{cor}{Corollary}
\newtheorem{lem}{Lemma}
\theoremstyle{definition}
\newtheorem*{defn}{Definition}

\newcommand{\Z}{\mathbb Z}
\newcommand{\R}{\mathbb R}
\newcommand{\C}{\mathbb C}
\newcommand{\A}{\mathcal A}
\newcommand{\X}{\underline X}

\newcommand{\D}{\mathcal{D}}
\newcommand{\J}[1]{J^{#1}(\R,\R)}
\newcommand{\dd}[1]{\frac{\partial}{\partial #1}}
\renewcommand{\d}{\,d}
\newcommand{\Ci}{C^{\infty}}
\newcommand{\g}{\mathfrak g}
\newcommand{\h}{\mathfrak h}
\newcommand{\AS}{\mathfrak S}

\newcommand{\sll}{\mathfrak{sl}}
\newcommand{\gl}{\mathfrak{gl}}
\renewcommand{\aa}{\mathfrak a}
\newcommand{\n}{\mathfrak n}
\newcommand{\m}{\mathfrak m}

\newcommand{\codim}{\operatorname{codim}}
\newcommand{\im}{\operatorname{im}}
\newenvironment{smt}{\left(\begin{smallmatrix}}{\end{smallmatrix}\right)}
\newcommand{\arctg}{\operatorname{arctg}}
\renewcommand{\sinh}{\operatorname{sh}}
\renewcommand{\cosh}{\operatorname{ch}}

\title{Contact Lie algebras of vector fields on the plane}
\authors{Boris M Doubrov\\ Boris P Komrakov}
\address{International Sophus Lie Centre\\ PO Box 70, 220123 Minsk, Belarus}
\email{Doubrov@islc.minsk.by\\Komrakov@islc.minsk.by}

\begin{document}
\begin{abstract}
  The paper is devoted to the complete classification of all real Lie
  algebras of contact vector fields on the first jet space of
  one-dimensional submanifolds in the plane. This completes Sophus Lie's
  classification of all possible Lie algebras of contact symmetries
  for ordinary differential equations. As a main tool we use the
  abstract theory of filtered and graded Lie algebras. We also
  describe all differential and integral invariants of new Lie
  algebras found in the paper and discuss the infinite-dimensional
  case.
\end{abstract}

\maketitlepage
\section{Introduction}

The problem of describing all finite-dimensional Lie algebras of vector
fields is a starting point for the symmetry analysis of ordinary differential
equations, because, having solved this problem, one finds all possible
algebras of contact symmetries for ordinary differential equations.

Over the complex numbers this classification was done at the end of the last
century by Sophus Lie~\cite{lie2}. He showed that, with three exceptions,
all Lie algebras of contact vector fields, viewed up to equivalence, are
lifts of Lie algebras of vector fields on the plane. The largest algebra
of the three exceptions (so-called {\em irreducible algebras of contact
vector fields\/}) is the algebra of contact symmetries of the equation
$y'''=0$ and is isomorphic to $\mathfrak{sp}(4,\C)$, while the other two
are its subalgebras of dimension 6 and~7.

In this paper we show that the problem of describing algebras of vector
fields can be formulated in a natural way in terms of filtered and graded Lie
algebras. This allows not only to give a new up-to-date proof of Sophus
Lie's classification, which is as yet missing in the literature, but also
to solve this problem over the field of real numbers. It turns out that
in the real case there are 8~irreducible contact Lie algebras of vector
fields on the plane, and one of them involves an arbitrary parameter.

Lie algebras of vector fields on the plane
were also classified (both in real and complex case) by Sophus
Lie~\cite{lie1}, so that the description of irreducible Lie algebras of
vector fields on the plane which is given in the present paper, basically
concludes the description of all finite-dimensional contact Lie algebras
of vector fields over the field of real numbers.

It should be noted that the problem of finding all irreducible contact
Lie algebras over the real numbers was also considered by F.~Engel
in~\cite{engel}, which is mentioned in Sophus Lie's three-volume
treatise~\cite{lie2} (volume 3, chapter~29, pages~760--761). P~Olver, in his
recently published book~\cite{olv}, cites this problem as unsolved.

\section{Jet space}

\subsection{Contact vector fields}
Let $M=\J1$ be the set of 1--jets of mappings from~$\R$ to~$\R$, and let~$\pi$
denote the natural projection $\J1\to\R^2$. We fix a coordinate system
$(x,y,z)$ on~$M$ in which the 1--jet of the mapping $f\co\R\to\R$ at the
point~$x_0$ has the coordinates $(x_0,f(x_0),f'(x_0))$. The projection~$\pi$
has in these coordinates the form $\pi\co (x,y,z)\mapsto (x,y)$.

We can introduce a natural {\em contact structure\/}~$M$. Indeed, there is
a two-dimensional distribution~$C$ on~$M$ which is not completely integrable
and has the property that all its integral curves whose projection onto the
plane is diffeomorphic, are precisely the curves of the form $(x,f(x),f'(x))$
with $f\in\Ci(I)$, $I\subset\R$. In terms of coordinates, this distribution
is given by the vector fields~$\dd z$ and $\dd x+z\dd y$ or, alternatively,
by the differential 1--form $\omega=\d y-z\d x$. A (local)
diffeomorphism~$\phi$ of the manifold~$M$ is said to be \emph{contact} if
$\phi$~preserves the contact distribution~$C$, ie, if
$d_p\phi(C_p)=C_{\phi(p)}$ for all $p\in M$. A vector field on~$M$ is called
contact if it generates a local one-parameter transformation group that
consists solely of contact diffeomorphisms. It is easy to show that a vector
field~$X$ is contact if and only if $L_X\omega=\lambda\omega$ for some
smooth function~$\lambda$.

If $X$ is a contact vector field, then the function $f=\omega(X)$ is called
the {\em characteristic function\/} of~$X$. It completely determines the
field~$X$, which in this case is denoted by~$X_f$ and has the form
$$X_f=-\frac{\partial f}{\partial z}\dd x + \left(f-z\frac{\partial
    f}{\partial z}\right)\dd y+ \left(\frac{\partial f}{\partial
    x}+z\frac{\partial f}{\partial y}\right)\dd z.$$

The mapping $f\mapsto X_f$ establishes an isomorphism between the space
of all smooth functions and that of contact vector fields on~$M$. This allows
to make the space $\Ci(M)$ into a Lie algebra by letting
$\{f,g\}=\omega([X_f,X_g])$.

\subsection{Prolongation operations} If $\phi$ is a (local) diffeomorphism
of the plane, then there exists a unique local contact transformation
$\phi^{(1)}\co\J1\to \J1$ such that the following diagram is commutative:
$$\begin{CD} \J1@>\phi^{(1)}>>\J1\\ @V\pi VV @V\pi VV \\
  \R^2@>\phi>>\R^2
\end{CD}$$
The transformation~$\phi^{(1)}$ is then called the (\emph{first})
{\em prolongation of the diffeomorphism\/}~$\phi$ and, in terms of
coordinates, has the from
$$\phi^{(1)}\co(x,y,z)\mapsto
\left(A(x,y),B(x,y),\frac{B_x+B_yz}{A_x+A_yz}\right).$$
Similarly, for any vector field~$X$ on the plane there exists a unique
contact vector field~$X^{(1)}$ on $\J1$ such that $\pi_*(X^{(1)})=X$. This
vector field~$X^{(1)}$ is called the (\emph{first}) {\em prolongation of the
vector field\/}~$X$ and has the form
$$X^{(1)}=A(x,y)\dd x + B(x,y)\dd y + (B_yz^2+(B_x-A_y)z-A_x)\dd z.$$
Its characteristic function is $B(x,y)-A(x,y)z$.

The mapping $X\mapsto X^{(1)}$ is an embedding of the Lie algebra of vector
field on the plane into the Lie algebra of contact vector fields on~$\J1$.
The contact vector fields that lie in the image of this mapping are called
\emph{point} contact vector fields. Point vector fields~$Y$ are characterized
by the following two equivalent properties:
\begin{enumerate}
\item any point vector field~$Y$ is an infinitesimal symmetry of
the \emph{vertical distribution\/}~$V$ on~$\J1$ ($V_p=\ker d_p\pi$);
\item the characteristic function of~$Y$ is linear in~$z$.
\end{enumerate}

\subsection{Reducible Lie algebras of contact vector fields}

\begin{defn} A Lie algebra~$\g$ of contact vector fields is called
  \emph{reducible} if there is a local contact diffeomorphism~$\phi$ such
  that the Lie algebra $\phi_*(\g)$ consists only of point vector fields.
  Otherwise, $\g$~is said to be \emph{irreducible}.
\end{defn}

\begin{thm}\label{th1}
  A Lie algebra~$\g$ of contact vector fields is irreducible if and only if
  it preserves no one-dimensional subdistribution of the contact distribution.
\end{thm}
\begin{proof}
  Every Lie algebra that consists of point vector fields preserves the
  vertical distribution~$V$, which is a one-dimensional subdistribution of the
  contact distribution~$C$. Consequently, any reducible Lie algebra of contact
  vector fields also preserves a one-dimensional subdistribution of~$C$.

  Conversely, let $\g$ be a Lie algebra of vector fields that preserves some
  one-dimensional subdistribution~$E$ of the contact distribution. If
  $A$ and~$B$ are two functionally independent first integrals of~$E$, then,
  as one can easily verify, the local diffeomorphism
  $$\phi\co (x,y,z)\mapsto
  \left(A,B,\frac{B_z}{A_z}=\frac{B_x+zB_y}{A_x+zA_y}\right)$$
  is contact and transforms the vertical distribution~$V$ to~$E$.
  It follows that the Lie algebra $\phi^{-1}_*(\g)$ preserves the vertical
  distribution and hence consists of point vector fields.
\end{proof}

\begin{cor} Any irreducible Lie algebra of contact vector fields is
  transitive at a point in general position.
\end{cor}
\begin{proof}
  Let $\g$ be an irreducible Lie algebra of contact vector fields. For an
  arbitrary point $p\in\J1$, we let $\g(p)=\{X_p\mid X\in\g\}$ and define
  $r=\max_{p\in\J1} \dim\g(p)$ and $U=\{p\in\J1\mid \dim\g(p)=r\}$. Then $U$
  is obviously an open subset in~$\J1$.

  The Lie algebra~$\g$ is transitive at point in general position if and only
  if $r=3$. Assume the contrary. Then the subspaces~$\g(p)$ form a completely
  integrable distribution~$E$ in~$U$ which is invariant under~$\g$. Consider
  the following two possibilities:

  1$^\circ$: $r=2$\qua 
  Then the intersection $E_p\cap C_p$ is one-dimensional at
  the points in general position, and this family of subspaces forms a
  one-dimensional subdistribution of the contact distribution which is
  invariant under~$\g$.

  2$^\circ$: $r<2$\qua In this case $E$~can be locally embedded into a
  two-dimensional completely integrable distribution which, as follows from
  its construction, is also invariant under~$\g$. Then, arguing as in the
  previous case, we conclude that the Lie algebra~$\g$ preserves a
  one-dimensional subdistribution of the contact distribution.
\end{proof}

In this paper we restrict ourselves to a local description of
finite-dimensional Lie algebras of contact vector fields at a point in
general position. In particular, from now on we shall assume that all
irreducible algebras of contact vector fields are transitive.

\section{An algebraic model of contact homogeneous space}

Let $\g$ be a transitive Lie algebra of contact vector fields on~$M=\J1$, let
$o$~be an arbitrary point in~$\J1$, and let $\g_0=\g_o$ be the subalgebra
of~$\g$ that consists of all vector fields in~$\g$ vanishing at the point~$o$.
It is easy to show that the subalgebra~$\g_0$ is \emph{effective}, ie,
contains no nonzero ideals of~$\g$ (see, for example,
\cite[Theorem 10.1]{herm}).

We can identify $T_oM$ with $\g/\g_0$ in the obvious way. Then $C_o$ is
identified with a certain submodule~$W$ of the $\g_0$--module $\g/\g_0$.
Since the distribution~$C$ is not completely integrable, the subspace
$\{x\in\g\mid x+\g_0\in W\}$ will not be closed with respect to the
multiplication in~$\g$.

We define a decreasing chain of subspace in~$\g$ as follows: $\g_p=\g$ for
all $p\le-2$,
$$\g_{-1}=\{x\in\g\mid x+\g_0\in W\},\quad \g_{p+1}=\{x\in\g_p\mid
[x,\g_{-1}]\subset \g_p\}\text{ for all }p\ge0.$$
It is easily shown by induction that $[\g_p,\g_q]\subset\g_{p+q}$ for all
$p,q\in\Z$, so that the family of subspaces $\{\g_p\}_{p\in\Z}$ defines a
filtration of the Lie algebra~$\g$.

\begin{defn} A filtered Lie algebra~$\g$ is called a \emph{contact Lie
  algebra} if
  \begin{enumerate}[a)]
  \item $\g_p=\g$ for all $p\le-2$;
  \item $\codim_{\g}\g_{-1}=1$, $\codim_{\g}\g_{0}=3$, and
    $[\g_{-1},\g_{-1}]+\g_{-1}=\g$;
  \item $\g_{p+1}=\{x\in\g_p\mid [x,\g_{-1}]\subset \g_p\}$ for all
    $p\ge0$;
  \item $\cap_{p\in\Z}\g_p=\{0\}$.
  \end{enumerate}
  Two contact Lie algebras are said to be isomorphic if they are isomorphic
  as filtered Lie algebras.
\end{defn}

Show that any transitive Lie algebra~$\g$ of contact vector fields is a
contact Lie algebra with respect to the above filtration. The properties
a) and c) follow immediately from the way that the filtration in~$\g$
is introduced. Let us prove b). Since the contact distribution
has codimension~1, we get $\codim_{\g}\g_{-1}=1$. From
transitivity of $\g$ on $M$ we get $\codim_{\g}\g_0=3$. Next,
since $C$ is not completely integrable, the subspace
$[\g_{-1},\g_{-1}]+\g_{-1}$ is strictly greater than
$\g_{-1}$ and, hence, is equal to $\g$.

Finally, Let $\aa=\cap_{p\in\Z}\g_p$. Then, obviously, $\aa$ is an ideal
in~$\g$ contained in~$\g_0$. Now since $\g_0$ is an effective subalgebra,
it follows that $\aa=\{0\}$. This proves~d).

Conversely, let $\g$ be an arbitrary finite-dimensional contact Lie algebra.
Then the pair $(\g,\g_0)$ determines a unique (up to local equivalence)
realization of~$\g$ as a transitive Lie algebra of vector fields on~$\R^3$.
And the subspace~$\g_{-1}$ allows us to define a $\g$--invariant
two-dimensional distribution on~$\R^3$ which is not completely integrable.
Therefore, the Lie algebra~$\g$ admits a unique (up to local equivalence)
realization as a transitive Lie algebra of contact vector fields.

Thus, the local classification of finite-dimensional transitive Lie algebras
of contact vector fields on~$\J1$ is equivalent to the classification (up to
isomorphism) of the corresponding contact Lie algebras. Observe that the
latter problem is algebraic and, as we shall see later, can be solved by
purely algebraic means.

All $\g$--invariant distributions on~$\J1$ are in one-to-one correspondence
with the submodules of the $\g_0$--module $\g/\g_0$. In particular, the
contact distribution corresponds to the submodule $\g_{-1}/\g_0$. From
Theorem~\ref{th1} now easily follows the next algebraic criterion for the
irreducibility of~$\g$.

\begin{lem}
\label{lem1}
  A transitive Lie algebra~$\g$ of contact vector fields is irreducible if
  and only if the $\g_0$--module $\g_{-1}/\g_0$ is irreducible.
\end{lem}

\section{Graded contact Lie algebras}

The major tool in the study of filtered Lie algebras is to consider graded
Lie algebras associated with them. As we shall see later on, with a few
exceptions, irreducible contact Lie algebras can be completely restored from
their associated graded Lie algebras.

\begin{defn}
  A $\Z$--graded Lie algebra $\h=\sum_{p\in\Z}\h_p$ is called a {\em graded
  contact Lie algebra\/} if
  \begin{enumerate}[a)]
  \item $\h_p=\{0\}$ for all $p<-2$;
  \item $\dim\h_{-1}=2$, $\dim\h_{-2}=1$, and $[\h_{-1},\h_{-1}]=\h_{-2}$;
  \item $\{x\in\h_p\mid [x,\h_{-1}]=0\}=\{0\}$ for all $p\ge0$.
  \end{enumerate}
\end{defn}

If $\g$ is a contact Lie algebra, then it is clear that the associated
graded Lie algebra $\h=\sum_p \g_p/\g_{p+1}$ satisfies all three conditions
in the above definition and is therefore a graded contact Lie algebra.
Moreover, if $\g$ is a filtered Lie algebra such that the associated graded
Lie algebra~$\h$ is contact and $\cap_p \g_p=\{0\}$, then it is easy to show
that $\g$ itself is a contact Lie algebra.

The concept of irreducibility for contact Lie algebras can be carried over to
graded contact Lie algebras. From Lemma~\ref{lem1} it immediately follows
that a contact Lie algebra~$\g$ is irreducible if and only if so is the
$\h_0$--module $\h_{-1}$ in the corresponding graded contact Lie algebra~$\h$.
The graded contact Lie algebras that satisfy this condition will be called
\emph{irreducible}.

The classification of all irreducible graded contact Lie algebras can be
carried out using the methods developed in the works of Tanaka~\cite{ta1,ta2}.
Slightly modifying the terminology of those papers (see also~\cite{stern}),
we introduce the concept of transitive graded Lie algebra, which
generalizes the concept of graded contact Lie algebra.

\begin{defn} A graded Lie algebra $\g=\oplus\g_p$ is said to be
\emph{transitive} if it satisfies the following conditions:
\begin{enumerate}[(i)]
\item there exists a natural number $\mu\in\mathbb N$ such that
$\g_{-p}=\{0\}$ for all $p>\mu$;
\item $[\g_{-1},\g_{-p}]=\g_{-p-1}$ for all $p\geq1$;
\item if $x\in\g_p$ for $p\geq0$ and $[x,\g_{-1}]=\{0\}$, then $x=0$.
\end{enumerate}
\end{defn}

It immediately follows from this definition that $\m=\bigoplus\limits_{p<0}
\g_p$ is the graded nilpotent Lie algebra generated by~$\g_{-1}$. Following
Tanaka~\cite{ta1,ta2}, we shall call graded nilpotent Lie algebras of this
kind \emph{fundamental}. In particular, the fundamental nilpotent Lie algebra
corresponding to a graded contact Lie algebra is none other than the
three-dimensional Heisenberg algebra.

Let $\m=\bigoplus\limits_{p<0}\g_p$ be an arbitrary fundamental graded Lie
algebra. Then, as was shown by Tanaka~\cite{ta1}, there exists a unique
transitive graded Lie algebra $\g(\m)=\bigoplus\limits_{p\in\mathbb Z}
\g_p(\m)$ that satisfies the following conditions:
\begin{enumerate}[(1)]
\item $\g_p(\m)=\g_p$ for $p<0$;
\item $\g(\m)$ is the largest among all transitive graded Lie
algebras satisfying condition~(1).
\end{enumerate}

This Lie algebra~$\g(\m)$ is called the (\emph{algebraic}) \emph{extension}
of~$\g$. In particular, any transitive graded Lie algebra~$\g$ may be
identified with a graded subalgebra of~$\g(\m)$, where
$\m=\bigoplus\limits_{p<0}\g_p$.

The Lie algebra~$\g(\m)$ has a clear geometrical meaning. Namely, let $M$~be
a connected Lie group with Lie algebra~$\m$, and let~$D$ be a left-invariant
distribution on~$M$ such that $D_e=\g_{-1}$. Denote by~$\A$ the Lie algebra
of all germs of infinitesimal symmetries of~$D$ at the identity element~$e$
of~$M$. Consider the following two subspaces in~$\A$:
\begin{align*}
\A_0=&\{\X\in\A\mid X_e=0\}\\
\A_{-1}=&\{\X\in\A\mid X_e\in D_e\}
\end{align*}
where $\X$ denotes the germ of the vector field~$X$ at the point~$e$. Now let
\begin{align*}
  A_{-p-1}=&[\A_{-p},\A_{-1}]\text{ for all }p\geq1\\
  A_p=&\{\xi\in\A_{p-1}\mid[\xi,\A_{-1}]\subset\A_{p-1}\}\text{ for all }
  p\geq1.
\end{align*}
Then the family of subspaces $\{\A_p\}_{p\in\mathbb Z}$ forms a decreasing
filtration of the Lie algebra~$\A$, and $\g(\m)$ can be identified with the
associated graded algebra, ie, $\g_p(\m)\equiv\A_p/\A_{p+1}$ for all
$p\in\mathbb Z$.

This geometrical interpretation allows to describe, without difficulty, the
structure of~$\g(\m)$ in the case that we are interested in, namely in the
case of graded contact Lie algebras.

Let $\n$ be the three-dimensional Heisenberg algebra and $\n_{-2}=[\n,\n]$,
while $\n_{-1}$ is a two-dimensional subspace complementary to~$[\n,\n]$.
In this case we may assume without loss of generality that $D$~is precisely
the contact distribution on~$\J1$.

Using the description of all infinitesimal symmetries of the contact
distribution, it is not hard to determine the structure of the Lie
algebra~$\g(\n)$. It can be identified with the space of polynomials in
$x,y,z$ with the bracket operation given by $X_{\{f,g\}}=[X_f,X_g]$. The
space~$\g_p(\n)$ consists of all homogeneous polynomials of degree~$p+2$,
assuming that the variables $x,y,z$ are of degree 1, 2 and~1 respectively.
For example,
\begin{align*}
\g_{-2}&(\n)=\langle1\rangle\\
\g_{-1}&(\n)=\langle x,z\rangle\\
\g_0&(\n)=\langle x^2,xz,z^2,y\rangle.
\end{align*}

We shall now fix some fundamental graded Lie algebra~$\m$ and describe how
one can classify all finite-dimensional graded subalgebras~$\h$ of the Lie
algebra~$\g(\m)$ such that $\h_{-p}=\m_{-p}$ for all $p\geq0$.

In what follows we shall always assume that $\h_{-p}=\m_{-p}=\g_{-p}(\m)$ for
all $p<0$. Suppose that for some $k\in\mathbb N\cup\{0\}$ we have a
collection of subspaces $\h_i\subset\g_i(\m)$, $i=0,\dots,k$, such that
$[\h_p,\h_q]\subset\h_{p+q}$ $\forall p,q\leq k$, $p+q\leq k$.
Using induction, we define a sequence of subspaces $\h_{k+1},\h_{k+2},\dots$
as follows:
$$\h_{p+1}=\{x\in\g_{p+1}(\m)\mid[x,\h_{-1}]\subset\h_p\}$$
for all $p\geq k$. It can be easily shown that
$$\g(\m,\h_0,\dots,\h_k)=\bigoplus_{p\in\mathbb Z}\h_p$$
is a graded subalgebra of~$\g(\m)$. This subalgebra is called the
\emph{extension} of the collection $(\h_0,\dots,\h_k)$. Note that
$\g(\m,\h_0, \dots,\h_k)$ is the largest of all graded subalgebras whose
$i$th grading space coincides with~$\h_i$ for all $i\leq k$.

One the other hand, we can associate $(\h_0,\dots,\h_k)$ with the graded
subalgebra $\tilde\g(\m,\h_0,\dots,\h_k)$ generated by these subspaces.

Now let $\h$ be an arbitrary graded subalgebra of~$\g(\m)$ such that
$\h_{-p}=\g_{-p}(\m)$ for all $p>0$. Then, obviously, for any $k\geq0$ we
have
$$\tilde\g(\m,\h_0,\dots,\h_k)\subset\h\subset\g(\m,\h_0,\dots,\h_k).$$
Finally, notice that $\g_0(\m)$ is precisely the algebra of all derivations
of~$\m$ that preserve the grading (see~\cite{ta1}), and all the subspaces
$\h_p\subset\g_p(\m)$ are invariant under the natural action of~$\h_0$
on~$\g_p(\m)$. Based on these remarks, the following algorithm for the
classification of the desired kind of subalgebras in~$\g(\m)$ suggests
itself.

\medskip
\noindent
\textbf{Step~I}\qua Describe, up to conjugation, all subalgebras
$\h_0\subset\g_0(\m)=\operatorname{Der}(\m)$. Go to Step~III.

\medskip
\noindent
\textbf{Step~II}\qua Suppose that for some $k\in\mathbb N\cup\{0\}$, a
collection of subspaces $\h_i\subset\g_i(\m)$, $i=0,\dots,k$, is already
constructed such that
\begin{align*}
{\rm(i)}&\quad[\h_p,\h_q]\subset\h_{p+q}\quad\forall p,q\leq k,\ p+q\leq k\\
{\rm(ii)}&\quad\dim\tilde\g(\m,\h_0,\dots,\h_k)<\infty.
\end{align*}
Let
\begin{align*}
&\tilde\g_{k+1}(\m,\h_0,\dots,\h_k)=\bigoplus_{\substack{i+j=k+1\\
1\leq i,j\leq k}}[\h_i,\h_j],\\
&\g_{k+1}(\m,\h_0,\dots,\h_k)=\{x\in\g_{k+1}(\m)\mid[x,\h_{-1}]\subset\h_k\}.
\end{align*}
At this point we describe all $\h_0$-invariant subspaces~$\h_{k+1}$ in
$\g_{k+1}(\m,\h_0,\dots,\h_k)$ such that
\begin{align*}
{\rm(i)}&\quad\tilde\g_{k+1}(\m,\h_0,\dots,\h_k)\subset\h_{k+1}\\
{\rm(ii)}&\quad\dim\tilde\g(\m,\h_0,\dots,\h_{k+1})<\infty.
\end{align*}

\smallskip
\noindent
\textbf{Step III}\qua Find the subalgebras $\tilde\g(\m,\h_0,\dots,\h_{k+1})$
and $\g(\m,\h_0,\dots,\h_{k+1})$. If these subalgebras are not the same, go to
Step~II. If, however, they coincide, then
$$\h=\tilde\g(\m,\h_0,\dots,\h_{k+1})=\g(\m,\h_0,\dots,\h_{k+1})$$
is one of the desired subalgebras.

\smallskip
Now we shall use this algorithm to classify all irreducible graded contact
Lie algebras over the field of real numbers.

\begin{thm}
\label{thm1}
  Let $\n$ denote the three-dimensional real Heisenberg algebra,
  considered as a graded Lie algebra, and let
  $\g(\n)$ be the universal extension of~$\n$. Then any finite-dimensional
  irreducible graded contact Lie algebra~$\h$ is isomorphic to one and only
  one of the following subalgebras of~$\g(\n)$:
  \begin{enumerate}[$1^\circ$]
  \item $\langle 1,x,y,z,x^2,xz,z^2,x(2y-xz),z(2y-xz),(2y-xz)^2\rangle$
  \item $\langle 1,x,y,z,x^2,xz,z^2\rangle$
  \item $\langle 1,x,z,x^2,xz,z^2\rangle$
  \item $\langle 1,x,z,x^2+z^2,2y-xz,x(x^2+z^2)-2z(2y-xz),
          z(x^2+z^2)+2x(2y-xz),(x^2+y^2)^2+4(2y-xz)^2\rangle$
  \item $\langle 1,x,z,x^2+z^2,2y-xz\rangle$
  \item $\langle 1,x,z,x^2+z^2+\alpha(2y-xz)\rangle$, $\alpha\geq0$
  \end{enumerate}
\end{thm}

\begin{proof}
Fix a basis $\{x,z\}$ in the space $\g_{-1}(\n)$. Then the action of
the elements of~$\g_0(\n)$ on $\g_{-1}(\n)$ is given by the following
matrices:
$$x^2\mapsto\begin{smt}0&2\\0&0\end{smt},\quad
xz\mapsto\begin{smt}-1&0\\0&1\end{smt},\quad
z^2\mapsto\begin{smt}0&0\\-2&0\end{smt},\quad
y\mapsto\begin{smt}-1&0\\0&0\end{smt}.$$
Therefore, the Lie algebra~$\g_0(\n)$ may be identified with $\gl(2,\R)$,
and the $\g_0(\n)$--module $\g_{-1}(\n)$ with the natural $\gl(2,\R)$--module.
\begin{lem}
\label{lem2}
  Any irreducible subalgebra of~$\gl(2,\R)$ is conjugate to one and only one
  of the following subalgebras:
  \begin{alignat*}{2}
  {\rm(i)}&\ \left\{
  \begin{pmatrix}\beta x&-x\\x&\beta x\end{pmatrix}
  \bigg|\ x\in\R\right\},\ \beta\geq0\quad &
  {\rm(ii)}&\ \left\{
  \begin{pmatrix}x&y\\-y&x\end{pmatrix}
  \bigg|\ x,y\in\R\right\}\\[2mm]
  {\rm(iii)}&\ \sll(2,\R)&{\rm(iv)}&\ \gl(2,\R)
  \end{alignat*}
\end{lem}
\begin{proof}
  If a subalgebra of $\gl(2,\R)$ is nonsolvable, then it is either
  three-di\-men\-si\-o\-nal and coincides with $\sll(2,\R)$, or
  four-dimensional and is equal to the whole of $\gl(2,\R)$. Any
  solvable irreducible subalgebra is commutative. If it is
  one-dimensional, then, as follows from the classification of real
  Jordan normal forms of $2\times2$ matrices, it is conjugate to the
  subalgebra~(i). If $\g$~is two-dimensional, it coincides with the
  centralizer of one of the Jordan normal forms, which implies that it
  is conjugate to the subalgebra~(ii).
\end{proof}

If we identify $\gl(2,\R)$ and $\g_0(\n)$, the subalgebras listed in
Lemma~\ref{lem2} are identified with the following subspaces
$\h_0\subset\g_0(\n)$:
\begin{alignat*}{2}
{\rm(i)\ }&\langle x^2+z^2+\alpha(2y-xz)\rangle,\ \alpha=2\beta\geq0\quad &
{\rm(ii)\ }&\langle x^2+z^2, 2y-xz\rangle\\[2mm]
{\rm(iii)\ }&\langle x^2, z^2, 2y-xz\rangle&
{\rm(iv)\ }&\langle x^2, xz, z^2, y\rangle
\end{alignat*}
Consider separately each one of these cases:

\smallskip\noindent
(i)\qua It is easily verified that in this case we have
$\g_1(\n,\h_0)=\{0\}$. Therefore $\h=\n\oplus\h_0$, and we arrive at the
algebra which is listed in the theorem under the number~$6^\circ$.

\smallskip\noindent
(ii)\qua Here we have
$$\g_1(\n,\h_0)=\langle x(x^2+z^2)-2z(2y-xz),
~z(x^2+z^2)+2x(2y-xz)\rangle,$$
and the action of the subalgebra~$\h_0$ on this space is irreducible.
Therefore, the space $\h_1\subset\g_1(\n,\h_0)$ is either zero or coincides
with the whole of~$\g_1(\n)$. In the former case we immediately find that
$\h=\n\oplus\h_0$ (subalgebra~$5^\circ$). In the second case the subalgebras
$\g(\n,\h_0,\h_1)$ and $\tilde\g(\n,\h_0,\h_1)$ coincide and are equal to
the subalgebra~$4^\circ$ of the theorem.

\smallskip\noindent
(iii)\qua Here $\g_1(\n,\h_0)=\langle x^3,x^2z,xz^2,z^3\rangle$, and the
$\h_0$--module $\g_1(\n,\h_0)$ is irreducible. Hence either we have
$\h_1=\{0\}$ and then $\h=\n\oplus\h_0$ (subalgebra~$3^\circ$), or
$\h_1=\g_1(\n,\h_0)$. In the latter case, however, the space~$\h_1$ generates
a finite-dimensional subalgebra.

\smallskip\noindent
(iv)\qua Here $\g_1(\n,\h_0)=\g_1(\n)$, and the $\h_0$--module $\g_1(\n)$ is a
sum of two irreducible submodules $W_1$ and~$W_2$ of the form
$$W_1=\langle x^3,x^2z,xz^2,z^3\rangle,\quad
W_2=\langle x(2y-xz),z(2y-xz)\rangle.$$
The submodule~$W_1$ generates a finite-dimensional subalgebra, so that
either $\h_1=\{0\}$ or $\h_1=W_2$. In the former case $\h=\n\oplus\h_0$
(subalgebra~$2^\circ$), while in the latter the subalgebras
$\g(\n,\h_0,\h_1)$ and $\tilde\g(\n,\h_0,\h_1)$ coincide and are equal to
the subalgebra~$1^\circ$ of the theorem.
\end{proof}

\section{Classification of contact Lie algebras}

In order to classify all finite-dimensional irreducible contact Lie algebras,
it will suffice to describe all filtered Lie algebras whose associated
graded Lie algebras are listed in Theorem~\ref{thm1}. To solve this latter
problem, we shall need the following result.

\begin{lem}\label{pas:l3}
  Let $\g$ be a finite-dimensional filtered Lie algebra, and $\h$ the
  associated graded Lie algebra. If there is an element $e\in\h_0$ such that
  $$[e,x_p]=px_p\quad\forall x_p\in\h_p$$
  then $\h$, viewed as a filtered Lie algebra, is isomorphic to~$\g$.
\end{lem}

\begin{proof} Suppose $e=\bar e+\g_1$ for some $\bar e\in\g_0$. For every
$p\in\mathbb Z$, consider the subspace
$$\g^p(\bar e)=\{x\in\g\mid[\bar e,x]=px\}.$$
It is easy to show that, $\g_p=\g^p(\bar e)\oplus\g_{p+1}$ for all
$p\in\mathbb Z$. Thus, the subspace $\g^p(\bar e)$ may be identified
with~$\h_p$, and since $[\g^i(\bar e),\g^j(\bar e)]\subset\g^{i+j}(\bar e)$,
this identification is in agreement with the structure of the Lie algebras
$\g$ and~$\h$. Hence, we have found an isomorphism of the Lie algebras
$\g$ and~$\h$ which is compatible with their filtrations.
\end{proof}

For the graded Lie algebras listed in Theorem~\ref{thm1} under the numbers
$1^\circ$, $2^\circ$, $4^\circ$,~$5^\circ$, we can choose~$e$ to be equal
to $xz-2y$, as this element is contained in all of these algebras. Then,
in view of~Lemma~\ref{pas:l3}, the description of the corresponding filtered
Lie algebras in these four cases is trivial. Consider the remaining two cases
$3^\circ$ and~$6^\circ$.

\smallskip\noindent
$3^\circ$\qua Let $\h$ be the graded Lie algebra that appears under the
number~$3^\circ$ in Theorem~\ref{thm1}, and let $\g$ be a contact Lie algebra
whose associated graded Lie algebra is isomorphic to~$\h$. Since
$\g_1=\{0\}$, the subalgebra~$\g_0$ can be identified with the
subalgebra~$\h_0$, which is isomorphic to $\sll(2,\R)$. Consider the
$\g_0$--module~$\g$. It is completely reducible, and its decomposition into a
sum of irreducible submodules has the form: $\g=V_{-2}\oplus
V_{-1}\oplus\g_0$, where the submodule~$V_{-2}$ is one-dimensional and is a
complement of~$\g_{-1}$, while the submodule~$V_{-1}$ is two-dimensional and
complements~$\g_0$ in~$\g_{-1}$. Therefore, the submodules $V_{-p}$, $p=1,2$
can be identified with the subspaces~$\h_{-p}$ of the graded Lie algebra~$\h$,
which allows to identify $\g$ and~$\h$ as vector spaces.

The structure of the Lie algebra~$\g$ is completely determined by the mappings
$\alpha\co V_{-2}\times V_{-1}\to\g$ and
$\beta\co V_{-1}\wedge V_{-1}\to\g$ defined as restrictions of the
bracket operation in~$\g$ to the corresponding subspaces. The Jacobi identity
shows that these mappings are both $\g_0$--invariant. Since the $\g_0$--module
$V_{-1}\wedge V_{-1}$ is one-dimensional and trivial, we have
$\im\beta\subset V_{-2}$. Similarly, $\im\alpha\subset V_{-1}$. Now,
computing the Jacobi identity for the basis vectors of $V_{-2}$
and~$V_{-1}$, we find that the mapping~$\alpha$ is zero. Thus, the
identification of the spaces $\g$ and~$\h$ is in agreement with the
Lie algebra structures of these spaces, so that the Lie algebra~$\g$ is
isomorphic to~$\h$, viewed as a filtered Lie algebra.

\smallskip\noindent
$6^\circ$\qua As in the above case, we can identify $\g_0$ and~$\h_0$. Now,
since the $\h_0$--modules $\h_0$ and $\h_{-1}$ are not isomorphic for any
value of~$\alpha$, we conclude that $\g_{-1}$ contains a $\g_0$--invariant
subspace~$V_{-1}$ which is a complement of~$\g_0$. Choose a basis~$\{e\}$
for~$\g_0$ and a basis $\{u_1,u_2\}$ for~$V_{-1}$ in such a way that
\begin{align*}
  [e,u_1]&=\alpha u_1-u_2\\
  [e,u_2]&=u_1+\alpha u_2.
\end{align*}
Then the elements $e,u_1,u_2$, together with the element $u_3=[u_1,u_2]$, will,
obviously, form a basis of~$\g$, and $[e,u_3]=2\alpha u_3$. Furthermore,
checking the Jacobi identity, we find that in case $\alpha\ne0$ we have
$[u_1,u_2]=[u_1,u_3]=0$, and the Lie algebra~$\g$ is isomorphic to~$\h$,
viewed as a filtered Lie algebra. If $\alpha=0$, we have $[u_1,u_2]=
\beta u_3$ and $[u_1,u_3]=-\beta u_2$ for some $\beta\in\R$. Note that the
parameters $\beta$ and~$x^2\beta$ with $x\in\R^*$ give here isomorphic Lie
algebras, whatever the value of~$x$ may be. Therefore, up to isomorphism of
contact Lie algebras we may assume that $\beta=0,\pm1$. If $\beta=0$, we
find that $\g$ is again isomorphic to~$\h$, viewed as a filtered Lie algebra.
If $\beta=1$ or $\beta=-1$, the Lie algebra~$\g$ is isomorphic to $\gl(2,\R)$
or~$\mathfrak u(2)$ respectively, while the subalgebras~$\g$ can be written,
under this identification, in matrix form as follows:
$$\begin{pmatrix} x & x \\ -x & x \end{pmatrix},\quad x\in\R.$$

Summing up what has been said, we obtain the following result:

\begin{thm}
  Any finite-dimensional irreducible contact Lie algebra is isomorphic to one
  and only one of the following:
  \begin{itemize}
  \item[\rm I] any of the graded contact Lie algebras listed in
  Theorem~\ref{thm1}, if they are viewed as filtered Lie algebras;

  \item[\rm II.1] $\g=\gl(2,\R)$, where $\g_p=\{0\}$ for $p\geq1$,
$$\g_0=\left\{\left.\begin{pmatrix} x & x \\ -x & x
  \end{pmatrix}\right|\,x\in\R\right\},\quad
\g_{-1}=\left\{\left.\begin{pmatrix} x+y & x+z \\ z-x & x-y
  \end{pmatrix}\right|\,x,y,z\in\R\right\};$$

  \item[\rm II.2] $\g=\mathfrak u(2)$, where $\g_p=\{0\}$ for $p\geq1$,
$$\g_0=\left\{\left.\begin{pmatrix} x & x \\ -x & x
  \end{pmatrix}\right|\,x\in\R\right\},\quad
\g_{-1}=\left\{\left.\begin{pmatrix} x+iy & x+iz \\ iz-x & x-iy
  \end{pmatrix}\right|\,x,y,z\in\R\right\}.$$
  \end{itemize}
\end{thm}

From now on, to refer to irreducible contact algebras of type~I, we shall
employ the notation I.$n$, where $n$ is the number of the corresponding
graded contact Lie algebra in Theorem~\ref{thm1}.

\section{Applications}

\subsection{} Now we shall find
explicit representations in contact vector
fields for the Lie algebras of vector fields described above.
Note that the mapping $f\mapsto X_f$ that maps an arbitrary function~$f$
of the variables $x,y,z$ into the vector field whose characteristic function
is~$f$, defines an embedding of the Lie algebra $\g(\n)$ into the algebra of
all contact vector fields. In this way we immediately obtain the explicit
representations in vector fields for those contact algebras~$\g$
which are isomorphic to their corresponding graded algebras.

Below we list three different representations of the space of characteristic
functions for each of the contact algebras II.1 and II.2:
\begin{itemize}
\item[II.1]
  \begin{itemize}
  \item[(a)] $\langle (2y-xz)^2+1, x-z(2y-xz), z+x(2y-xz), x^2+z^2
    \rangle$
  \item[(b)] $\langle x^2+z^2, 2x(2y-xz)+z(x^2+z^2+4),
    2z(2y-xz)-x(x^2+z^2+4), 16+4(2y-xz)^2+(x^2+z^2)^2 \rangle$
  \item[(c)] $\langle 1,z,\sqrt{1+z^2}\sinh x, \sqrt{1+z^2}\cosh x \rangle$
  \end{itemize}
\item[II.2]
  \begin{itemize}
  \item[(a)] $\langle (2y-xz)^2+1, x+z(2y-xz), z-x(2y-xz), x^2+z^2
    \rangle$
  \item[(b)] $\langle x^2+z^2, 2x(2y-xz)+z(x^2+z^2-4),
    2z(2y-xz)-x(x^2+z^2-4), 16+4(2y-xz)^2+(x^2+z^2)^2 \rangle$
  \item[(c)] $\langle 1,z,\sqrt{1-z^2}\sin x, \sqrt{1-z^2}\cos x \rangle$
  \end{itemize}
\end{itemize}

In particular, from the representations (a) and~(b) it follows that
these two algebras of contact vector fields can both be embedded into the
10--dimensional algebra~I.1 and into the 8--dimensional algebra~I.4. The
representations~(c) are notable for the fact that the characteristic
functions here are independent of~$y$.

\subsection{} Consider the set of all contact vector fields of the form~$X_f$,
where the function~$f$ has the form $f=ay+g(x,z)$ with $a\in\R$ and
$g$~being an arbitrary function of $x,z$. It is easy to show that this
condition is equivalent to the requirement that $X_f$~be an infinitesimal
symmetry of the one-dimensional distribution~$E$ generated by the vector
field $\dd{y}$. Thus we see that this space of vector fields forms an
infinite-dimensional subalgebra~$\AS$ of the Lie algebra of all contact
vector fields.

Consider the projection $\pi\co\AS\to\D(\R^2)$ given by
$$\quad X_f=-g_z\dd x +(g-zg_z)\dd y+(g_x+az)\dd z\mapsto -g_z\dd x +
(g_x+az)\dd z.$$
It is easily verified that this mapping is a homomorphism of Lie algebras
whose kernel is one-dimensional and is generated by~$X_1$, while its image
coincides with the set of all vector fields on the plane that preserve,
up to a constant factor, the volume form $\omega=dx\wedge dz$ on the plane:
\begin{equation}\label{eq:pr}
\pi(\AS)=\left\{X\in\D(\R^2)\mid L_x(\omega)=\lambda\omega,\
  \lambda\in\R\right\}.
\end{equation}

Thus, with every Lie algebra of contact vector fields that preserves a
one-dimensional distribution complementary to the contact one, we can
associate a Lie algebra of vector fields on the plane. Conversely, the
inverse image of any subalgebra of the Lie algebra~(\ref{eq:pr}) of vector
fields on the plane is some Lie algebra of contact vector fields in the jet
space.

Note that all irreducible Lie algebras of contact vector fields, except
I.1 and~I.4, preserve a one-dimensional distribution complementary to the
contact distribution, and hence can be embedded into~$\AS$. The corresponding
Lie algebras of vector fields on the plane are as follows:
\begin{itemize}
\item[I.2] the Lie algebra corresponding to the group of affine
  transformations of the plane;
\item[I.3] the Lie algebra corresponding to the group of equi-affine
  transformations of the plane (ie, area-preserving affine transformations);
\item[I.5] the Lie algebra corresponding to the group of similitude
  transformations;
\item[I.6] $\langle \dd x, \dd z, (\beta x-z)\dd x+(x+\beta z)\dd
  z\rangle$, $\beta=\alpha/2$ (if $\alpha=0$, this Lie algebra corresponds to
  the group of Euclidean transformations);
\item[II.1] the Lie algebra corresponding to the group of all
  transformations of the hyperbolic plane;
\item[II.2] the Lie algebra corresponding to the group of rotations of the
  sphere.
\end{itemize}

\subsection{} The above correspondence allows to describe without any
difficulty all differential and integral invariants for all Lie algebras~$\g$
of contact vector fields that satisfy the following conditions:
\begin{itemize}
\item[(a)] $\g$ preserves a one-dimensional distribution complementary to
the contact one;
\item[(b)] $\g\ni X_1=\dd y$.
\end{itemize}

Indeed, let $(x,y_0=y,y_1=z,y_2,\dots,y_n)$ be the standard coordinate system
in the space $J^n(\R^2)$ of $n$th jets of curves on the
plane. (See~\cite{olv} for definition of jet spaces and notions of
differential and integral invariants.) Denote by $\g^{(n)}$ the $n$th
prolongation of the Lie algebra $\g$.  It then follows
from the condition~(b) that the differential and integral invariants of~$\g$
that have the order~$n$ are independent of~$y$ and may be considered on the
manifold of the trajectories of the vector field~$X_1^{(n)}$. These
trajectories are given by the equations $y=\operatorname{const}$ and can be
parametrized by the coordinates $(x,y_1,\dots,y_n)$. Furthermore, it turns
out that if $n\ge2$, the action of the algebra~$\g^{(n)}$ on that quotient
manifold is equivalent to the action of the Lie algebra $\pi(\g)^{(n-1)}$ on
the space of $(n-1)$th jets, and the mapping $J^n(\R^2)\to J^{(n-1)}(\R^2)$
that establishes this equivalence has the form:
$$(x,y_1,\dots,y_n)\mapsto(x,y_0,\dots,y_{n-1}).$$
Therefore, all differential and integral invariants of~$\g$ may be derived
from the invariants of~$\pi(\g)$ by substituting $y_{i+1}$ instead of~$y_i$
for $i\ge0$.

We remark that Sophus Lie~\cite{lie1} found all invariants for those Lie
algebras of vector fields on the plane that correspond to the cases I.2
and~I.3. The invariants of the 10--dimensional irreducible Lie algebra~I.1
were computed in~\cite{olv} over the complex numbers, and they remain
unchanged on passing to the real case. Now we shall specify nontrivial
integral and differential invariants of the least order for the rest of
irreducible contact Lie algebras of vector fields; all other invariants can
be derived from these by means of differentiation (see~\cite{olv}).

$$
\begin{array}{|c|c|c|}
\hline
 & \text{Differential invariant} & \text{Integral invariant}\\
\hline
\text{I.4} & \frac{P}{Q^{8/3}} & \frac{Q^{1/3}dx}{y_2^2+1} \\[2mm]
\text{I.5} & \frac{(1+y_2^2)y_4-3y_2y_3^2}{y_2^2} &
\frac{y_3dx}{1+y_2^2} \\[2mm]
\text{I.6} & \frac{y_3e^{-\beta\arctg y_2}}{(1+y_2^2)^{3/2}} &
e^{\beta\arctg y_2}(1+y_2)^{1/2}dx \\[2mm]
\text{II.1} &
\frac{\sqrt{1+y_1^2}((1+y_1^2)y_3-3y_1y_2^2-y_1(1+y_1^2)^2)}%
{((1+y_1^2)^2+y_2^2)^{3/2}} &
\left(\frac{(1+y_1^2)^2+y_2^2}{1+y_1^2}\right)^{1/2}dx \\[2mm]
\text{II.2} &
\frac{\sqrt{1-y_1^2}((1-y_1^2)y_3+3y_1y_2^2+y_1(1-y_1^2)^2)}%
{((1-y_1^2)^2+y_2^2)^{3/2}} &
\left(\frac{(1-y_1^2)^2+y_2^2}{1-y_1^2}\right)^{1/2}dx\\[2mm]
\hline
\end{array}
$$
where
\begin{align*}
P&=(y_2^2+1)^2\big(QD^2(Q)-\tfrac{7}{6}D(Q)\big)+2(y_2^2+1)y_2y_3QD(Q)-\\
&\,\hskip5cm-\big(9(y_2^2+1)y_2y_4-\tfrac{1}{2}(9y_2^2-19)y_3^2\big)Q^2\\
Q&=9(y_2^2+1)^2y_5-90(y_2^2+1)y_2y_3y_4+5(27y_2^2-5)y_3^3\\
D&=\dd x +y_1\dd y +\dots + y_7 \dd {y_6}\quad\text{(the operator of
total differentiation).}
\end{align*}
For the algebras II.1 and~II.2, we have chosen here their representations
in contact vector fields that appear earlier under the letter~(c).

Notice that all contact Lie algebras listed in the table above
are reducible over the
field of complex numbers.  Hence, for each of these algebras there
exists a certain complex analytic contact transformation which takes
it to one of the known canonical forms for contact Lie algebras over
$\C$.  Thus, the inverse thatsformation (prolonged as many times as
needed) brings known invariants to the invariants of the initial
Lie algebra.

For example, the contact transformation
$$T\co (x,y,z)\co(x,y,z)\mapsto(x+iz,-2iy+1/2(x^2+2ixz+z^2),x-iz)$$
takes the contact Lie algebra I.5 to the algebra with the following space
(over $\C$) of characteristic functions:
$$\langle 1, x, y, z, xz\rangle.$$
This  contact Lie algebra is reducible and
is the first prolongation of the following
Lie algebra of vector fields on the plain:
$$\left\langle \dd{x}, x\dd{x}, y\dd{y}, \dd{y}, x\dd{y}\right\rangle.$$
The differential invariants of the least order for this Lie algebra
were computed already (see, for example,~\cite{olv}) and have the
form:
\begin{center}
\begin{tabular}{cc}
differential invariant: & $\frac{y_2y_4}{y_3^2}$ \\
integral invariant: & $\frac{y_3}{y_2}\,dx$.
\end{tabular}
\end{center}
The third prolongation of the inverse transformation $T^{-1}$ takes
these invariants to those given in the table above. In the similar way
we can compute invariants for other contact Lie algebras of vector
fields given in the table.

\subsection{} Consider the problem of classifying those infinite-dimensional
subalgebras in the Lie algebra of contact vector fields that correspond to
Lie pseudo-groups of contact transformations (ie, those that can be defined
with the help of a finite number of differential equations; see~\cite{stern}).
As in the finite-dimensional case, all these algebras can be naturally
divided into two classes: reducible ones, which are actually extensions of
infinite-dimensional Lie algebras of vector fields on the plane, and
irreducible ones. Over the field of complex numbers all irreducible
infinite-dimensional Lie algebras of contact vector fields were described by
Sophus Lie~\cite{lie3}, who showed that, apart from the Lie algebra of all
contact vector fields, there exist exactly two infinite-dimensional
irreducible subalgebras, namely, the Lie algebra $\AS=\{X_{ay+g(x,z)}\}$,
which we already mentioned earlier, and its commutant
$[\AS,\AS]=\{X_{g(x,z)}\}$.

The methods for the description of contact Lie algebras that have been
developed in this paper, can be easily generalized to the
infinite-dimensional case. In particular, with these Lie algebras we can
again associate graded contact Lie algebras that can be embedded into the
universal extension~$\g(\n)$ of the three-dimensional Heisenberg algebra.
Let $\h$ be an infinite-dimensional graded subalgebra in~$\g(\n)$ such that
$\h_p=\g_p(\n)$ for $p<0$ and such that the $\h_0$--module $\h_{-1}$ is
irreducible. Then all possible types of subalgebras~$\h_0$ over the real
numbers are listed in Lemma~\ref{lem2}. As follows from the proof of
Theorem~\ref{thm1}, in the cases (i) and~(ii) the Lie algebra~$\h$ is
finite-dimensional. The consideration of the remaining cases (iii) and~(iv)
is the same over the complex and real numbers, and gives the
infinite-dimensional Lie algebras of contact vector fields described above.
Thus the classification of irreducible infinite-dimensional Lie algebras of
contact vector fields remains unchanged on passing from the complex to the
real case.

\end{document}